\documentclass[letterpaper, 10 pt, conference]{ieeeconf}  % Comment this line out if you need a4paper

\usepackage{amsmath,amssymb}
\usepackage{graphicx} 
\usepackage{epsfig}
\usepackage{amsfonts}
\usepackage{amstext}

\newcommand{\done}{\hspace*{\fill} $\Box$}
\newcommand{\RE}{I\!\!R}
\newcommand{\IN}{I\!\!N}

\newtheorem{theorem}{Theorem}
\newtheorem{remark}{Remark}

\newtheorem{definition}{Definition}

\newtheorem{example}{Example}

\begin{document}

\author{Luis Rodrigues\\
Department of Electrical and Computer Engineering\\
 Concordia University, 1515 St. Catherine Street, Montr\'eal, QC H3G 2W1, Canada\\}

\title{\LARGE \bf Inverse Optimal Control with Discount Factor for Continuous and Discrete-Time Control-Affine Systems and Reinforcement Learning}
\maketitle

\begin{abstract}
This paper addresses the inverse optimal control problem of finding the state weighting function that leads to a quadratic value function when the cost on the input is fixed to be quadratic.
The paper focuses on a class of infinite horizon discrete-time and continuous-time optimal control problems whose dynamics are control-affine and whose cost is quadratic in the input.
The optimal control policy for this problem is the projection of minus the gradient of the value function onto the space formed by all feasible control directions.
This projection points along the control direction of steepest decrease of the value function.
For discrete-time systems and a quadratic value function the optimal control law can be obtained as the solution of a regularized least squares program, which corresponds to a receding horizon control with a single step ahead.
For the single input case and a quadratic value function the solution for small weights in the control energy is interpreted as a control policy that at each step brings the trajectories of the system as close as possible to the origin, as measured by an appropriate norm. 
Conditions under which the optimal control law is linear are also stated.
Additionally, the paper offers a mapping of the optimal control formulation to an equivalent reinforcement learning formulation.
%For the case of linear dynamics and quadratic cost the paper proposes a model-free reinforcement l%earning control strategy using semidefinite programming.
Examples show the application of the theoretical results.
\end{abstract}

%\begin{keywords}
%optimal control; Lie point symmetries; Hamilton-Jacobi-Bellman equation; Lyapunov function;
%\end{keywords}

\section{Introduction}
Unconstrained linear optimal control problems have been completely solved and have shown its importance in several applications.
However, the open literature in nonlinear optimal control problems is much sparser than the one for linear optimal control.
For a detailed literature survey on optimal control see \cite{Athans1966,Bryson1969,Sussman1997,Liberzon2012} and references therein.
Optimal control problems lead to the Bellman equation in discrete-time and to the Hamilton-Jacobi-Bellman equation in continuous-time.
It is often difficult to obtain solutions of these equations for nonlinear systems, which has motivated the early study of inverse optimal control laws for continuous-time systems \cite{Kalman1964,Thau1967,Yokoyama1972,Anderson1973}.
One class of nonlinear problems that have been studied at greater depth is the class of infinite horizon optimal control for dynamical systems that are control-affine and whose cost is quadratic in the input \cite{Popescu2005}.
Under certain conditions these systems are feedback linearizable, which was the focus of reference \cite{Freeman1996}.
The work in \cite{Freeman1996} used a control Lyapunov function obtained for the linearized dynamics and yielded an analytical solution for a stabilizing controller.
The existing work in optimal and inverse optimal control can be divided into three classes:
\begin{enumerate}
\item the use of a cost function without a discount factor (e.g, \cite{Lee2021,Prasanna2019, Dierks2010,Ornelas2010}), 
\item imposing the constraint that the weighting function on the state must be quadratic (e.g.,\cite{Popescu2005}),
\item techniques that limit the number of states and inputs (e.g., \cite{Rodrigues2016,Clelland2013}).
\end{enumerate}
It is important to note that none of the references \cite{Thau1967,Yokoyama1972,Anderson1973,Popescu2005,Freeman1996,Lee2021,Prasanna2019, Dierks2010,Ornelas2010,Rodrigues2016,Clelland2013} have solved the infinite horizon inverse optimal control problem for both discrete-time and continuous-time systems. They also do not provide any link to the area of reinforcement learning. A description of the connections of reinforcement learning and optimal control, including an example in linear quadratic control, can be found in \cite{Lewis2009}. However, inverse optimal control was not addressed in \cite{Lewis2009}.

The focus of this paper is on a unified approach to inverse optimal control of both discrete-time and continuous-time control-affine systems, which can also be helpful to parameterize the value function in a reinforcement learning strategy.
\iffalse
Although it has been shown before that control-affine systems with quadratic control cost have an analytical solution for the optimal control law, to the best of the author's knowledge this paper offers a novel geometric interpretation for these solutions.
\fi

This paper addresses inverse optimal control  for both continuous-time and discrete-time control-affine dynamics with a discount factor with the following contributions:
\begin{itemize}
\item there is no structural constraint in the parameterization of the state weighting function in the three formal results presented in this paper,
\item the results are valid for an arbitrary number of states and inputs and for a cost function with a discount factor,
\item a geometric interpretation is offered for the optimal control policy based on the notion of inner product (see figure \ref{interpretation}),
\item an interpretation of the optimal control law as the solution of a regularized least squares problem is offered for discrete-time optimal control with a quadratic value function,
\item the inverse optimal control results are connected to the formulation of reinforcement learning using an exponential mapping.
\end{itemize}

The paper is organized as follows. Section \ref{ocproblem} addresses the discrete-time optimal control problem followed by inverse optimal discrete-time control in section \ref{inverseocp_discrete}. Then section \ref{continuousoptimal} addresses the continuous-time optimal control problem followed by inverse optimal continuous-time control in section \ref{inverseocp_continuous}.
After section \ref{reinforcementlearning} addressing the connections of the proposed approach to reinforcement learning, the examples and conclusions close the paper.

\section{Discrete-Time System Model and Optimal Control Problem}\label{ocproblem}
\noindent We consider a discrete-time nonlinear dynamical system that is affine in the control input written as
\begin{eqnarray}\label{system}
x(k+1) &=& f(x(k)) + g(x(k)) u(k)
\end{eqnarray}
where $x\in\RE^n, u\in\RE^m$, $f(x)$ is Lipschitz,  $g(x)$ has maximum rank for all $x\neq0\in\RE^n$, and $f(0)=0$.
We define a control policy as a set of control inputs of the form $\pi_{u_k}=\{u(k)~:k\in\IN_0\}$, where $\IN_0$ is the set of non-negative integers.
It is assumed that the system (\ref{system}) is stabilizable according to the following definition.
\vspace{10pt}

\begin{definition}
The system (\ref{system}) is stabilizable if there exists a control policy that can asymptotically take the trajectories of the system to the origin for any initial condition $x(0)$.
\end{definition}
\vspace{10pt}

\begin{figure}[t] 
\centerline{ \resizebox{120mm}{!}{\includegraphics{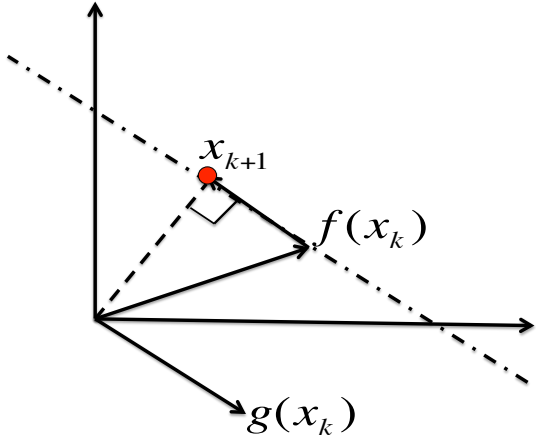}} }
\caption{Geometric interpretation of optimal control for a quadratic value function and a single input when $R(x(k))<<\gamma \|g(x(k))\|_P^2$}
\label{interpretation}
\end{figure}

The cost of a control policy $\pi_{u_k}$ when the system starts at the state $x(k)$ will be measured by the function
\begin{equation}\label{cost}
V(x(k),\pi_{u_k})=\sum_{i=k}^{\infty}\gamma^{i-k}\left[Q(x(i))+u^T(i)R(x(i))u(i)\right]
\end{equation}
where $R(x(i))=R(x(i))^T> 0$ and $Q(0)=0$.
Conditions under which $Q(x(i))\ge 0$ for all $x(i)\in\RE^n$ will be stated in Theorem \ref{inverseoptimality}.
The parameter $0<\gamma\le 1$ is a discount factor.
Splitting the sum in (\ref{cost}) into the term for $i=k$ plus an infinite sum starting at $i=k+1$, the cost (\ref{cost}) becomes
\begin{eqnarray}\label{Bellman}
V(x(k),\pi_{u_k})=Q(x(k))+u^T(k)R(x(k))u(k)+\nonumber\\
+\gamma V(x(k+1),\pi_{u_{k+1}}).
\end{eqnarray}
Equation (\ref{Bellman}) is known as the Bellman equation.
Let us define the value function as $V^*(x(k))=\inf_{\pi_{u_k}} V(x(k),\pi_{u_k})$, where the minimization is constrained to the dynamics in equation (\ref{system}).
Bellman's principle of optimality states that a sub-trajectory of an optimal trajectory whose initial state is $x(0)$, and that starts at a state $x(k)$, is itself an optimal trajectory with $x(k)$ as the new initial condition.
Therefore, according to Bellman's principle of optimality and (\ref{Bellman}),
\begin{eqnarray}\label{Bellmanoptimality}
V^*(x(k))&=&Q(x(k))+\nonumber\\
&+&\inf_{u(k)}\left[u^T(k)R(x(k))u(k)+\gamma V^*(x(k+1))\right]\nonumber\\
\end{eqnarray}
where $x(k+1)$ is determined by equation (\ref{system}).
Notice that the term $Q(x(k))$ was removed from the minimization in equation (\ref{Bellmanoptimality}) because it does not depend on the control input $u(k)$.
From equation (\ref{Bellmanoptimality}) we see that the optimal control policy is the recursive solution of the optimization problem
\begin{equation}
\begin{array}{rclr}
\inf_{u(k)} &\left[u^T(k)R(x(k))u(k)+\gamma V^*(x(k+1))\right]&\\
 s.t. &x(k+1) = f(x(k)) + g(x(k)) u(k)\label{ocp}&\\
\end{array}
\end{equation}
for all current and future values of time $k$.
This problem can be rewritten in the form 
\begin{equation}
\begin{array}{rclr}
\inf_{u(k)} &\left[u^T(k)R(x(k))u(k)+\gamma V^*(x(k+1))\right]&\\
 s.t. &f(x(k))=- g(x(k)) u(k) + x(k+1)& \label{newocp}\\
\end{array}
\end{equation}
where $x(k+1)$ can be interpreted as the error of modelling $f(x(k))$ as the function $-g(x(k))u(k)$, which is linear in $u(k)$.
The necessary condition of optimality is
\begin{equation}\label{necessarycondition}
2R(x(k))u(k)+\gamma g^T(x(k))\nabla {V^*}\left(x(k+1)\right)=0,
\end{equation}
where equation (\ref{system}) can be used to replace $x(k+1)$.
Solving for the control input $u(k)$ in (\ref{necessarycondition}) yields
\begin{equation}\label{discreteu}
u(k)=-\frac{\gamma}{2}R^{-1}(x(k))g^T(x(k))\nabla {V^*}\left(x(k+1)\right),
\end{equation}
which is a projected gradient descent control law.
In other words, the term $g(x(k)))u(k)$ in (\ref{system}) with $u(k)$ from (\ref{discreteu}) is a scaled projection of the gradient of $V^*(x(k+1))$ onto $g(x(k))$, which gives the coordinates of a vector along the space of columns of $g(x(k)))$ that leads to decreasing values of the value function $V^*(x(k+1))$.
The control law (\ref{discreteu}) is also proportional to minus the directional derivative of $V^*(x(k+1))$ along the direction of $g(x(k))$ when there is a single input and $R$ is constant, 
The value function $V^*(x)$ is the solution of equation (\ref{Bellmanoptimality}) with $u(k)$ replaced by the solution of (\ref{discreteu}).
If the optimal value function is quadratic then it admits the parameterization
\begin{equation}\label{quadraticvalue}
V^*(x(k+1))=x^T(k+1)Px(k+1)
\end{equation}
where $P=P^T>0$.
In this case the optimization (\ref{newocp}) is a regularized least squares problem, which
is equivalent to a single step ahead receding horizon optimal control problem.
The analytic solution of equation (\ref{discreteu}) then becomes 
\begin{eqnarray}\label{quadraticsolution}
u(k)=-\gamma\left[R(x(k))+\gamma g^T(x(k))Pg(x(k))\right]^{-1}\cdot\nonumber\\
\cdot g^T(x(k))Pf(x(k)).
\end{eqnarray}
For the case of a single input, i.e., $m=1$, defining the inner product
\begin{equation}\label{innerproduct}
\left<f(x),g(x)\right>_P=g^T(x)Pf(x)=f^T(x)Pg(x)
\end{equation}
and the norm
\begin{equation}\label{norm}
\|g(x)\|_P=\sqrt{\left<g(x),g(x)\right>_P},
\end{equation}
the solution (\ref{quadraticsolution}) can be rewritten as
\begin{eqnarray}\label{newquadraticsolution}
u(k)=-\gamma\left[R(x(k))+\gamma \|g(x(k))\|_P^2\right]^{-1}\cdot\nonumber\\
\cdot\left<f(x(k)),g(x(k))\right>_P,
\end{eqnarray}
which can be cast in the form
\begin{equation}\label{u}
u(k)=-\left[1+\frac{R(x(k))}{\gamma \|g(x(k))\|_P^2}\right]^{-1}\left<f(x(k)),\frac{g(x(k))}{\|g(x(k))\|_P^2}\right>_P
\end{equation}
Therefore, when $m=1$ and $R(x(k))<<\gamma \|g(x(k))\|_P^2$ we can interpret the solution (\ref{u}) as the negative value of the projection of the drift vector $f(x(k))$ onto the control direction $g(x(k))/\|g(x(k))\|_P$.
Note that in this case the control (\ref{u}) is is the input that brings the state $x(k+1)$ as close to zero as possible along the control direction (see figure \ref{interpretation}), where the distance is measured by the norm (\ref{norm}).
To the best of the author's knowledge this is a new interpretation of the optimal control policy (\ref{newquadraticsolution}).
Note that when $P$ is the identity matrix one can implement the control law (\ref{u}) with only knowledge of the mathematical model and without the need of a parameterizing matrix.
The closed-loop system dynamics for a single input is obtained by replacing (\ref{u}) into (\ref{system}) yielding
\begin{eqnarray}\label{closedloop}
x(k+1)=f(x(k))-g(x(k))\left[1+\frac{R(x(k))}{\gamma \|g(x(k))\|_P^2}\right]^{-1}\cdot\nonumber\\
\cdot\left<f(x(k)),\frac{g(x(k))}{\|g(x(k))\|_P^2}\right>_P.
\end{eqnarray}
\vspace{10pt}
\begin{remark}\label{remark1}
Notice that the optimal control input (\ref{u}) is guaranteed to stabilize the closed-loop system when $\gamma=1$ and $Q(x)$ is positive definite because in this case the value function $V^*(x)$  is a Lyapunov function according to the Bellman equation (\ref{Bellmanoptimality}).
\end{remark}
\vspace{10pt}
\begin{remark}\label{remark_geometry}
The geometric interpretation in figure \ref{interpretation} can be extended to the case of multiple inputs. The most straightforward extension is to the case when the matrix $g(x)$ has orthogonal columns. In that case, the inner product (\ref{innerproduct}) can be applied to each column $g_i(x)$. The term $g^T(x(k))Pf(x(k))$ then leads to the coordinates of the projection of $f(x)$ onto the orthogonal basis formed by the columns of $g(x)$, whose span represents the space of all possible control directions. 
\end{remark}

\section{Discrete-Time Inverse Optimal Control}\label{inverseocp_discrete}
\noindent The choice of the function $Q(x(k))$ for a particular problem is the task of the control designer.
In the interest of using the analytic solution (\ref{quadraticsolution}), and to be able to claim its optimality, the following question must be answered: which choice (if any) of $Q(x(k)),~Q(0)=0,$ would lead to a solution of the Bellman equation (\ref{Bellmanoptimality}) that is a quadratic value function $V^*(x(k))=x(k)^TPx(k)$?
Notice that for such $Q(x(k))$ one would not only have an analytical solution for the control input, but also a clear geometric interpretation  (see figure \ref{interpretation}).
These are two important motivating points to investigate the answer to the posed question.
Choosing a specific parameterization of the function $V^*(x(k))$ and then solving the Bellman equation (\ref{Bellmanoptimality}) for $Q(x(k))$ is called the inverse optimal control problem.
\vspace{10pt}

\begin{theorem}\label{inverseoptimality}
If the function $Q(x(k))$ is chosen as
\begin{eqnarray}\label{qfinal}
Q(x)=x(k)^TPx(k)-\gamma f^T(x(k))Pf(x(k))+\nonumber\\
\gamma^2f^T(x(k))Pg(x(k))M^{-1}(x(k))g^T(x(k))Pf(x(k)),
\end{eqnarray}
where
\begin{equation}\label{M}
M(x(k))=R(x(k))+\gamma g^T(x(k))Pg(x(k)),
\end{equation}
\iffalse
\begin{equation}\label{q}
Q(x(k))=x(k)^TPx(k)-\gamma f^T(x(k))P\left[f(x(k))+g(x(k))u(k)\right],
\end{equation}
for $P=P^T>0$ with $u(k)$ given by (\ref{quadraticsolution}) 
\fi
with $P=P^T>0$ and $u(k)$ given by (\ref{quadraticsolution}) , then $V^*(x(k))=x(k)^TPx(k)$ satisfies the Bellman equation (\ref{Bellmanoptimality}).
Furthermore, for the case of a single input, the function $Q(x(k))$ in (\ref{qfinal}) can be rewritten as
\begin{equation}\label{qsolution}
Q(x(k))=\|x(k)\|_P^2-\gamma Q_2(x(k)),
\end{equation}
where
\begin{eqnarray}\label{distancesquared}
Q_2(x(k))=\|f(x(k))\|_P^2-\left[1+\frac{R(x(k))}{\gamma \|g(x(k))\|_P^2}\right]^{-1}\cdot\nonumber\\
\cdot\left|\left<f(x(k)),\frac{g(x(k))}{\|g(x(k))\|_P}\right>_P\right|^2,
\end{eqnarray}
with  inner product and norm as defined in (\ref{innerproduct}) and (\ref{norm}), respectively.
For $R(x(k))<<\gamma \|g(x(k))\|_P^2$ the expression (\ref{qsolution}) can be approximated as
\begin{eqnarray}\label{caserzero}
Q(x(k))&\approx&\|x(k)\|_P^2-\gamma\|x(k+1)\|_P^2\nonumber\\
&=&V^*(x(k))-\gamma V^*(x(k+1)).
\end{eqnarray}
Additionally, if $f(x)$ is Lipschitz continuous in $x$ with constant $L>0$  using the norm (\ref{norm}), then for  $\gamma\le L^{-1}$ the function $Q(x(k))$ in (\ref{qfinal}) is non-negative for all $x(t)\in\RE^n$ and $Q(0)=0$.
\end{theorem}
\vspace{10pt}
\proof
Rearranging equation (\ref{quadraticsolution}) and multiplying by $u^T(k)$ on the left one can write
\begin{eqnarray}\label{newu}
u^T(k)\left[R(x(k))+\gamma g^T(x(k))Pg(x(k))\right]u(k)\nonumber\\
=-\gamma
u^T(k)g^T(x(k))Pf(x(k)).
\end{eqnarray}
Using equations (\ref{system}) and (\ref{newu}) and $Q(x(k))$ given by (\ref{qfinal}) one can see that $V^*(x(k))=x(k)^TPx(k)$ satisfies the equation 
(\ref{Bellmanoptimality}), which upon these substitutions becomes
\begin{eqnarray}\label{newoptimality}
x(k)^TPx(k)&=&Q(x(k))+\nonumber\\
\gamma f^T(x(k))Pf(x(k))&+&\gamma f^T(x(k))Pg(x(k))u(k)~
\end{eqnarray}
Replacing (\ref{quadraticsolution}) into (\ref{newoptimality}) yields (\ref{qfinal}).
For the case of a single input, replacing (\ref{u}) in  (\ref{qfinal}) yields (\ref{qsolution}).
When $R(x(k))<<\gamma \|g(x(k))\|_P^2$ the term $Q_2(x(k))$ in equation (\ref{qsolution}) can be approximated by a difference of squares. By the theorem of Pythagoras (see figure \ref{interpretation} for a depiction of the right triangle) this difference of squares is equal to the square of the distance of $x(k+1)$ to the origin, where the distance is measured by $\|x\|_P$ using the norm (\ref{norm}).
Therefore, one obtains the expression (\ref{caserzero}).
Finally, to prove under what conditions $Q(x(k))$ is non-negative we can see from equation (\ref{qsolution}) that the worst case for $Q(x(k))$ to potentially become negative would be when $R(x(k))\to\infty$.
For this case,
\begin{equation}\label{qsolutionworstcase}
Q(x(k))=\|x(k)\|_P^2-\gamma\|f(x(k))\|_P^2.
\end{equation}
Since $f$ is Lipschitz continuous, denoting the positive Lipschitz constant by $L$ we have
\begin{equation}\label{Lipschitz}
\|f(x(k))\|_P^2\le L\|x(k)\|_P^2,
\end{equation}
which when replaced in equation (\ref{qsolutionworstcase}) yields
\begin{equation}
Q(x(k))\ge\left(1-\gamma L\right)\|x(k)\|_P^2.
\end{equation}
Therefore, $Q(x(k))$ will be non-negative if $\gamma\le L^{-1}$.
Additionally, $Q(0)=0$ since $f(0)=0$.
\done
\vspace{10pt}
\begin{remark}
Theorem \ref{inverseoptimality} extends the results of reference \cite{Ornelas2010} to the case of a cost function with a discount factor $\gamma\neq 1$.
Additionally, theorem \ref{inverseoptimality} provides a condition on the Lipschitz constant of the function $f$ to guarantee that $Q(x(k))$ is non-negative. Furthermore, theorem \ref{inverseoptimality} does not need any condition of detectability as in reference \cite{Ornelas2010}.
\end{remark}
\vspace{10pt}
\begin{remark}
If $\gamma=1$ and $R(x(k))<<\gamma \|g(x(k))\|_P^2$ can be neglected, then when the approximate expression (\ref{caserzero}) is introduced into the summation (\ref{cost}) it leads to telescoping series. This series then yields a cost that penalizes the distance of the initial state to the origin provided that the optimal control is stabilizing and $\lim_{k\to\infty}x(k)\to 0$. 
\end{remark}
\vspace{10pt}

\begin{remark}
When $L\le 1$ then $\gamma=1$ is a possible choice that yields a guarantee of stability, as stated in remark \ref{remark1}.
\end{remark}

\section{Continuous-Time System Model and Optimal Control Problem}\label{ocproblem_continuous}\label{continuousoptimal}
\noindent We now consider a continuous-time nonlinear dynamical system that is affine in the control input written as
\begin{eqnarray}\label{continuous_system}
\dot x(t) &=& f(x(t)) + g(x(t))u(t)
\end{eqnarray}
where $x\in\RE^n, u\in\RE^m$, $f(x)$ is Lipschitz,  $g(x)$ has maximum rank  for all $x\neq0\in\RE^n$, and $f(0)=0$.
As before, we assume that the system (\ref{continuous_system}) is stabilizable.
The control policy $\pi_u$ is defined as $\pi_u=\{u(t)~:t\in[t_0,\infty)\}$.
The cost of a control policy $\pi_u$ when the system starts at the state $x(t_0)$ will be measured by the function
\begin{eqnarray}\label{continuous_cost}
V(x(t_0),\pi_u)=\nonumber\\
\int_{t_0}^{\infty}e^{-\gamma(t-t_0)}\left[Q(x(t))+u^T(t)R(x(t))u(t)\right]dt
\end{eqnarray}
where $R(x(t))=R^T(x(t))>0$ and $Q(0)=0$.
Conditions under which $Q(x(t))\ge 0$ for all $x\in\RE^n$ will be stated in Theorem \ref{inverseoptimality_continuous}.
The parameter $\gamma\ge 0$ is a discount factor.
The admissible control policies $\pi_u^*$ are all Lebesgue measurable functions $u(t),~t\in[t_0,\infty)$, for which the integral in equation (\ref{continuous_cost}) is finite.
Splitting the integral in (\ref{continuous_cost}) in two intervals $[t_0,t_0+T]$ and $[t_0+T,\infty)$, where $T>0$, and minimizing for $u\in\pi_u^*$ yields
\begin{eqnarray}\label{newcontinuouscost}
V^*(x(t_0))=\nonumber\\
\inf_{u\in\pi_u^*}\left[\int_{t_0}^{t_0+T}
e^{-\gamma(t-t_0)}J(t)dt+e^{-\gamma T}V^*(x(t_0+T))\right],
\end{eqnarray}
where
\begin{equation}\label{Jcost}
J(t)=\left[Q(x(t))+u^T(t)R(x(t))u(t)\right]
\end{equation}
Bellman's principle of optimality can now be applied in the limiting case when $T\to 0$.
To do that we perform a first order Taylor series approximation of the right hand side of equation (\ref{newcontinuouscost}), divide by $T$, and take the limit when $T\to 0$ (assuming the limit exists) to yield the Hamilton-Jacobi-Bellman equation
\begin{eqnarray}\label{HJB}
\gamma V^*(x(t))=Q(x(t))+~~~~~~~~~~~~~~~~~~~\nonumber\\
\inf_{u(t)}\left[u^T(t)R(x(t))u(t)+\nabla{V^*}^T\dot x(t)\right],
\end{eqnarray}
where $Q(x(t))$ was removed from the minimization because it does not depend on $u(t)$ and $\dot x(t)$ is given by (\ref{continuous_system}).
%\left(f(x(t))+g(x(t))u(t)\right)
The necessary condition of optimality for the control input yields the projected gradient descent law
\begin{equation}\label{necessaryconditionHJB}
u(t)=-\frac{1}{2}R^{-1}(x(t))g^T(x(t))\nabla V^*(x(t)).
\end{equation}
When equation (\ref{necessaryconditionHJB}) is compared to the discrete-time law (\ref{discreteu}) we observe that an important difference is that for (\ref{necessaryconditionHJB}) the control input at time $t$ only depends on the gradient of the value function $V^*(x(t))$ at the same time $t$.
For the single input case if $R(x(t))=0.5\|g((t)x)\|_P^2$ then the optimal control vector field becomes
\begin{equation}\label{singleinputu}
u(t)g(x(t))=-\left<\nabla V^*(x(t)),\frac{g(x(t))}{\|g(x(t))\|_P^2}\right>_Pg(x(t)),
\end{equation}
which is minus the projection of $\nabla V^*(x(t))$ onto the control vector $g(x(t))$.
This indicates the direction along $g(x(t))$ that leads to decreasing values of the value function $V^*(x(t))$.
The function $V^*(x(t))$ is the solution of equation (\ref{HJB}) with $u(t)$ replaced by (\ref{necessaryconditionHJB}).
If the optimal value function is quadratic then
\begin{equation}\label{quadraticvalueHJB}
V^*(x(t))=x^T(t)Px(t)
\end{equation}
where $P=P^T>0$.
For this case, equation (\ref{necessaryconditionHJB}) becomes
\begin{equation}\label{quadraticsolutionHJB}
u(t)=-R^{-1}(x(t))g^T(x(t))Px(t).
\end{equation}
For a single input, i.e., $m=1$, using the inner product (\ref{innerproduct}) and the norm (\ref{norm}) the solution (\ref{quadraticsolutionHJB}) can be rewritten as
\begin{equation}\label{newquadraticsolutionHJB}
u(t)=-R^{-1}(x(t))\left<g(x(t)),x(t)\right>_P,
\end{equation}
The closed-loop system dynamics for a single input is obtained by replacing (\ref{newquadraticsolutionHJB}) into (\ref{continuous_system}) yielding
\begin{equation}\label{closedloop_continuous}
\dot x(t)=f(x(t))-g(x(t))R^{-1}(x(t))\left<g(x(t)),x(t)\right>_P.
\end{equation}
Notice that the vector $x(t)$ points in the radial direction away from the origin.
Therefore, the control input (\ref{newquadraticsolutionHJB}) forces the closed-loop vector field $\dot x(t)$ in (\ref{closedloop_continuous}) to point  toward the closest point to the origin along the line with equation $f(x)+g(x)u$, where the distance is measured by the norm (\ref{norm}).
To the best of the author's knowledge this is a new interpretation of the optimal control solution (\ref{newquadraticsolutionHJB}).
\vspace{10pt}
\begin{remark}
As in the discrete-time case, it can be shown that the geometric interpretation of the optimal control law (\ref{newquadraticsolutionHJB}) can be extended to the case of multiple inputs (see remark \ref{remark_geometry}).
\end{remark}
\vspace{10pt}
\begin{remark}\label{remark3}
Notice that the optimal control input (\ref{newquadraticsolutionHJB}) is guaranteed to stabilize the closed-loop system when $\gamma=0$ and $Q(x(t))$ is positive definite because in this case the value function $V^*(x(t))$  is a Lyapunov function according to the Hamilton-Jacobi-Bellman equation (\ref{HJB}).
\end{remark}

\section{Continuous-Time Inverse Optimal Control}\label{inverseocp_continuous}
\noindent We now address the inverse optimality question: for what choice of $Q(x(t))$ is the solution of the Hamilton-Jacobi-Bellman equation (\ref{HJB}) a quadratic value function $V^*(x(t))=x(t)^TPx(t)$?
\vspace{10pt}

\begin{theorem}\label{inverseoptimality_continuous}
If the function $Q(x(t))$ is chosen as
\begin{eqnarray}\label{qcontinuous}
Q(x(t))=\gamma x(t)^TPx(t)-2x(t)^TPf(x(t))+\nonumber\\
 x(t)^TPg(x(t))R^{-1}(x(t))g^T(x(t))Px(t),
\end{eqnarray}
for $P=P^T>0$ with $u(t)$ given by (\ref{quadraticsolutionHJB}) then $V^*(x(t))=x(t)^TPx(t)$ satisfies the Hamilton-Jacobi-Bellman equation (\ref{HJB}).
Furthermore, if $R(x(t))$ is chosen as a constant matrix and $g(x(t))$ is also a constant matrix then the optimal control (\ref{quadraticsolutionHJB}) is linear.
For the case of a single input, the function $Q(x(t))$ in (\ref{qcontinuous}) can be rewritten as
\begin{eqnarray}\label{qsolution_continuous}
Q(x(t))=\gamma\|x(t)\|_P^2-2\left<f(x(t)),x(t)\right>_P+\nonumber\\
+R^{-1}(x(t))|\left<g(x(t)),x(t)\right>_P|^2.
\end{eqnarray}
Additionally, if $f(x)$ is Lipschitz continuous in $x$ with constant $L>0$ using the norm (\ref{norm}), then for  $\gamma\ge 2L$ the function $Q(x(t))$ in (\ref{qsolution_continuous}) is non-negative for all $x(t)\in\RE^n$ and $Q(0)=0$.
\end{theorem}
\vspace{10pt}
\proof
Replacing $u(t)$ by (\ref{quadraticsolutionHJB}) and $Q(x(t))$ by (\ref{qcontinuous}) one can see that $V^*(x(t))=x(t)^TPx(t)$ satisfies the HJB equation (\ref{HJB}).
The proof of the statement about the linearity of the control when both $R$ and $g$ are constant matrices is trivial from equation (\ref{quadraticsolutionHJB}).
For the single input case, using the notions of inner product (\ref{innerproduct}) and norm (\ref{norm}), the equation (\ref{qcontinuous}) can be rewritten as (\ref{qsolution_continuous}).
Since $f$ is Lipschitz continuous, denoting the positive Lipschitz constant by $L$ one can write
\begin{equation}\label{Lipschitz_continuous}
\|f(x(t))\|_P^2\le L\|x(t)\|_P^2.
\end{equation}
From equation (\ref{qsolution_continuous}) the worst case scenario for $Q(x)$ to potentially become negative is when $R(x)\to\infty$.
For this case, using the Lipschitz relation (\ref{Lipschitz_continuous}) and Schwartz's inequality 
\begin{equation}\label{Schwartz}
\left<f(x),x\right>_P\le\|f(x)\|_P\|x\|_P,
\end{equation}
yields
\begin{equation}
Q(x(t))\ge\left(\gamma-2L\right)\|x(t)\|_P^2.
\end{equation}
Therefore, $Q(x(t))$ will be non-negative if $\gamma\ge 2L$.
Additionally, $Q(0)=0$.
\done
\vspace{10pt}

\noindent As per remark \ref{remark3} it is convenient to choose $\gamma=0$ for stability purposes.
However, if $\gamma=0$ the condition $\gamma\ge 2L$ will not be satisfied because $L>0$.
The following theorem gives a condition under which $Q(x(t))$ is non-negative when $\gamma=0$ for the single input case.
\vspace{10pt}
\begin{theorem}\label{gamma0}
Assume the stabilizable dynamical system (\ref{continuous_system}) with a single input.
Let $\mathcal S_{g^\perp}$ be defined as
\begin{equation}\label{gperpset}
\mathcal S_{g^\perp}=\left\{x\in\RE^n~:~\left<g(x),x\right>_P=0\right\}.
\end{equation}
When $\gamma=0$ the function $Q(x)$ in equation (\ref{qsolution_continuous}) will be non-negative for all $x\in\RE^n$ if $\left<f(x),x\right>_P\le 0$ for all $x\in\mathcal S_{g^\perp}$. Furthermore, if $\left<f(x),x\right>_P< 0$ for all $x\in\mathcal S_{g^\perp},~x\neq 0$, then the  norm $\|x\|_P$ as defined in (\ref{norm}) will be decreasing with time if $R(x)$ is chosen such that
\begin{equation}\label{R}
\frac{\left<f(x),x\right>_P}{\left|\left<g(x),x\right>_P\right|^2}-R(x)>0,~\forall~x\in\RE^n.
\end{equation}
\end{theorem}
\vspace{10pt}
\proof
The proof that $Q(x)$ is non-negative is trivial from expression (\ref{qsolution_continuous}).
Proving that there is $R(x(t))$ small enough that the norm $\|x(t)\|_P$ is decreasing with time is equivalent to proving the same for the square of the norm, i.e., $\left<x(t),x(t)\right>_P$.
The time rate of change of the square of the norm is
\begin{equation}\label{squarednormrate}
z(t)=\frac{d\|x(t)\|_P^2}{dt}=2\left<x(t),\dot x(t)\right>_P.
\end{equation}
Replacing (\ref{closedloop_continuous}) into (\ref{squarednormrate}) yields
\begin{eqnarray}
%\frac{d\|x(t)\|_P^2}{dt}
z(t)=\nonumber
2\left[\left<f(x(t)),x(t)\right>_P-R^{-1}(x(t))\left|\left<g(x(t)),x(t)\right>_P\right|^2\right]
\end{eqnarray}
For any finite $R(x(t))$ when $x(t)\in\mathcal S_{g^\perp},~x\neq 0$ we conclude that $\|x(t)\|_P^2$ has negative derivative from the assumptions of the theorem.
When $x\notin\mathcal S_{g^\perp}$, then $\|x(t)\|_P^2$ has negative derivative if $R(x(t))$ is chosen according to (\ref{R}).
\done
\vspace{10pt}
\begin{remark}
The result on the norm $\|x\|_P$ can be extended to the case of inputs bounded by an hypersphere if the number of inputs is the same as the number of states (see chapter 10 of reference \cite{Athans1966} for details).
\end{remark}

\section{Connections to Reinforcement Learning}\label{reinforcementlearning}
This paper considers the minimization of a cost function, which is the typical approach in optimal control.
Notice however that a maximization of a reward function would yield the same results if the reward was written as the exponential $r(x,\pi_u)=e^{-V(x,\pi_u)}$, where $V(x,\pi_u)$ denotes the expression (\ref{cost}) for discrete-time systems and the expression (\ref{continuous_cost}) for continuous-time systems. 
The function $r(x,\pi_u)$ can instead be interpreted as a probability of reward or probability of success given that the cost $V(x,\pi_u)$ is always non-negative, which leads to $r(x,\pi_u)\in(0,1]$.
The reward (or utility) approach is the one used in reinforcement learning and also in optimal decision making in finance and economics.
Therefore, the results of this paper can also be applied to these areas using the exponential mapping outlined above.
Future work in reinforcement learning can explore this connection to inverse optimal control and vice-versa.

\section{Examples}
The following examples show the application of the theoretical results developed in the paper.
\vspace{10pt}
\begin{example}
For the second order system described by (\ref{system}) with state $x(k)=[x_1(k)~~x_2(k)]^T$,
\begin{eqnarray*}
f(x(k)) = \left[
\begin{array}{c}
-x_2(k)\sin(x_2(k))\\
-x_1(k)\cos(x_2(k))\sin(x_1(k))
\end{array}
\right],\\
g(x(k)) = \left[
\begin{array}{c}
0\\
1
\end{array}
\right],
\end{eqnarray*}
$P=I$, and $R(x(k))=R$, from therorem \ref{inverseoptimality} there exists a $Q(x(k))$ for which  $V^*(x(k))=\|x(k)\|^2$ is the optimal value funciton.
The value function in this case measures the distance of the state to the origin.
We now compute the optimal control $u(k)$ and the function $Q(x(k))$ (when $R\to 0$).
Notice that
\begin{eqnarray*}
\left|\left<f(x(k)),\frac{g(x(k))}{\|g(x(k))\|_P}\right>_P\right|^2=~~~~~~~~~~~~~~~~~~~~~~~~\nonumber\\
x_1^2(k)\cos^2(x_2(k))\sin^2(x_1(k)),\\
\\
\|f(x(k))\|^2=~~~~~~~~~~~~~~~~~~~~~~~~\nonumber\\
x_2^2(k)\sin^2(x_2(k))+x_1^2(k)\cos^2(x_2(k))\sin^2(x_1(k)),
\end{eqnarray*}
and the Lipschitz constant is $L=1$ because
\begin{eqnarray*}
\|f(x(k))\|^2\le\|x(k)\|^2.
\end{eqnarray*}
From theorem \ref{inverseoptimality} and remark \ref{remark1}, the optimal control input for the choice $\gamma=1$ is stabilizing and given by (\ref{u}), i.e.,
\begin{eqnarray*}
u(x(k))=\frac{x_1(k)\cos(x_2(k))\sin(x_1(k))}{1+R}.
\end{eqnarray*}
When $R\to 0$ this is a deadbeat control and the state-weighting function becomes
\begin{eqnarray*}
Q(x(k)) = x_1^2(k)+x_2^2(k)\left[1-\sin^2(x_2(k))\right]\ge 0.
\end{eqnarray*}
\end{example}
\vspace{10pt}
\begin{example}
For the second order system described by (\ref{continuous_system}) with state $x(t)=[x_1(t)~~x_2(t)]^T$,
\begin{eqnarray*}
f(x(t)) = \left[
\begin{array}{c}
x_2^3(t)-x_1(t)\\
-x_1(t)x_2^2(t)
\end{array}
\right],\quad
g(x(t)) = \left[
\begin{array}{c}
0\\
1
\end{array}
\right],
\end{eqnarray*}
$P=I$, and $R(x(t))=R$, from theorem \ref{inverseoptimality_continuous} there exists a $Q(x(t))$ for which  $V^*(x(t))=\|x(t)\|^2$ is the optimal value function.
The value function thus measures the distance of the state to the origin.
Notice that the functions $f$ and $g$ are such that they satisfy the conditions of theorem \ref{gamma0} because $\left<f(x),x\right>=-x_1^2\le 0$.
Therefore, one can choose $\gamma=0$ and obtain a non-negative function $Q(x)$.
According to remark \ref{remark3}, the choice of $\gamma=0$ guarantees that the optimal controller is stabilizing.
We now compute the optimal control $u(t)$ and the function $Q(x(t))$.
The optimal control input given by (\ref{quadraticsolutionHJB}) is a linear state feedback and is written as
\begin{eqnarray*}
u(x(t))=-R^{-1}x_2(t).
\end{eqnarray*}
From (\ref{qsolution_continuous}) the function $Q(x(t))$ for $\gamma=0$ is
\begin{eqnarray*}
Q(x(t)) = 2x_1^2(t)+R^{-1}x_2^2(t),
\end{eqnarray*}
which is a positive definite function.
\end{example}

\section{Conclusions}
\iffalse
This paper offered a new geometric interpretation of the analytical optimal control solution for systems that are affine in the control with a cost that is quadratic in the input. This new interpretation opens the door to future development of data-driven algorithms that can approximate the optimal solution in real-time. Furthermore, 
\fi
This paper presented an analytical solution of the inverse optimal control problem when the value function is quadratic. Both discrete-time and continuous-time systems were considered. Examples showed the application of the theory.

\bibliographystyle{IEEEtran}
\bibliography{ocpfmsbibliography}
\end{document}